\documentclass[11pt]{article}
\usepackage{amsfonts}
\usepackage{color}
\usepackage{amsmath,amssymb,comment}
\newtheorem{theo}{Theorem}[section]
\newtheorem{lem}[theo]{Lemma}
\newtheorem{cor}[theo]{Corollary}
\newtheorem{prop}[theo]{Proposition}

\newcommand{\mysection}[1]{\section{#1} \setcounter{equation}{0}}
\newcommand{\proof}{{\sc Proof.} \quad}
\newcommand{\proofc}{{\sc Proof} \ }
\newcommand{\be}{\begin{equation} \label}
\newcommand{\ee}{\end{equation}}
\newcommand{\bea}{\begin{eqnarray}\label}
\newcommand{\eea}{\end{eqnarray}}
\newcommand{\bas}{\begin{eqnarray*}}
\newcommand{\eas}{\end{eqnarray*}}
\newcommand{\bit}{\begin{itemize}}
\newcommand{\eit}{\end{itemize}}
\newcommand{\qed}{\hfill$\Box$ \vskip.2cm}
\newcommand{\nn}{\nonumber}
\newcommand{\R}{\mathbb{R}}

\newcommand{\pO}{\partial\Omega}

\newcommand{\hra}{\hookrightarrow}
\newcommand{\io}{\int_\Omega}
\newcommand{\na}{\nabla}

\newcommand{\lam}{\lambda}

\newcommand{\bom}{\overline{\Omega}}
\newcommand{\Om}{\Omega}

\newcommand{\wh}{\widehat}

\newcommand{\vs}{\vspace*}
\newcommand{\hs}{\hspace*}
\newcommand{\proj}{{\mathcal{P}}}

\newcommand{\abs}{\\[5pt]}

\newcommand{\tm}{T_{max}}

\newcommand{\Kf}{K_F}
\newcommand{\kf}{k_F}
\begin{document}
\enlargethispage{10mm}
\title{Conditional estimates in three-dimensional chemotaxis-Stokes systems\\
and application to a Keller-Segel-fluid model\\
accounting for gradient-dependent flux limitation
}
\author{
Michael Winkler\footnote{michael.winkler@math.uni-paderborn.de}\\
{\small Institut f\"ur Mathematik, Universit\"at Paderborn,  }\\
{\small Warburger Str.~100, 33098 Paderborn, Germany}
}
\date{}
\maketitle
\begin{abstract}
\noindent
  This manuscript is concerned with the Keller-Segel-Stokes system
  \bas
	\left\{ \begin{array}{lcl}
	n_t + u\cdot\na n &=& \Delta n - \na \cdot \Big( n F(|\na c|^2) \na c\Big), \\[1mm]
	c_t + u\cdot\na c &=& \Delta c - c + n, \\[1mm]
	u_t &=& \Delta u + \na P + n\na\Phi, 
	\qquad \na \cdot u=0,
	\end{array} \right.
	\qquad \qquad (\star)
  \eas 
  under no-flux/no-flux/Dirichlet boundary conditions in smoothly bounded three-dimensional domains,
  with given suitably regular functions $F$ and $\Phi$. Here in accordance with recent developments in the literature
  on refined modeling of chemotactic migration, the introduction of suitably decaying $F$ is supposed to adequately account for
  saturation mechanisms that limit cross-diffusive fluxes near regions of large signal gradients.
  In the context of such nonlinearities which suitably generalize the prototype given by $F(\xi)=\Kf (1+\xi)^{-\frac{\alpha}{2}}$,
  $\xi\ge 0$, with $\Kf>0$, 
  known results addressing a fluid-free parabolic-elliptic simplification of ($\star$) have identified the value 
  $\alpha_c=\frac{1}{2}$ as critical with regard to the occurrence of blow-up in the sense that some exploding solutions
  can be found when $\alpha<\frac{1}{2}$, whereas all suitably regular initial data give rise to global bounded solutions
  when $\alpha>\frac{1}{2}$.\abs
  The intention of the present study consists in making sure that 
  the latter feature of blow-up prevention by suitably strong flux limitation persists
  also in the more complex framework of the fully coupled chemotaxis-fluid system ($\star$).
  To achieve this,
  as a secondary objective of possibly independent interest the manuscript separately establishes 
  some conditional bounds for corresponding fluid fields and taxis gradients in a fairly general setting that
  particularly includes the subsystem of ($\star$) concerned with the evolution of $(c,u,P)$.
  These estimates relate respective regularity features to certain integrability properties of associated forcing terms,
  as in the context of ($\star$) essentially represented by the quantity $n$.\abs
  The application of this tool to the specific problem under consideration thereafter facilitates 
  the derivation of a result on global existence of bounded classical solutions to ($\star$) for widely arbitrary initial data
  actually within the entire range $\alpha>\frac{1}{2}$, and by means of an argument which appears to be signficantly 
  condensed when compared to reasonings pursued in previous works concerned with related problems.\abs
\end{abstract}
\newpage
\mysection{Introduction}\label{intro}
{\bf Flux-limited Keller-Segel systems in liquid environments.} \quad
Due to their relevance in refined models for chemotactic movement (\cite{bellomo_flim}, \cite{perthame},
\cite{bianchi_painter_sherratt}, \cite{bianchi_painter_sherratt2016}), dependencies of cell
migration rates on gradients of the unknown quantities have received increasing interest in the recent literature on
mathematical analysis of taxis-type cross-diffusion systems.
Exemplary studies are concerned with global solvability in models involving $p$-Laplace type cell diffusion operators
(\cite{bendahmane_et_al}, \cite{li_yan2020}, \cite{zhuang_wang_zheng_sining2020},
\cite{tao_weirun_li_yuxiang2019}, \cite{tao_weirun_li_yuxiang2020}, \cite{liu}), 
quite precisely identify critical parameters for the occurrence of blow-up in systems
simultaneously accounting for flux-limited cross-diffusion and diffusion mechanisms paralleling those in relativistic heat
equations (\cite{bellomo_win1}, \cite{bellomo_win2}, \cite{chiyoda_mizukami_yokota}, \cite{mizukami_ono_yokota}), 
or also discuss wave-like solution behavior and pattern-supporting instability features in various particular modeling
contexts based on flux-limiting mechanisms (\cite{arias_campos_soler}, \cite{calvez_perthame_yasuda}, \cite{perthame}).\abs
With regard to 
issues related to singular structure formation,
among the model classes most comprehensively understood seems to be the family of parabolic-elliptic systems given by
\be{PE}
	\left\{ \begin{array}{l}
	n_t = \Delta n - \na \cdot \Big( nF(|\na c|^2)\na c\Big), \\[1mm]
	0 = \Delta c - c + n,
	\end{array} \right.
\ee
where the given function $F$ describes gradient-dependent limitiation of crodd-diffusive
fluxes by suitably generalizing the prototype
\be{fp}
	F(\xi)=\Kf (1+\xi)^{-\frac{\alpha}{2}},
	\qquad \xi \ge 0,
\ee
with $\Kf>0$ and $\alpha>0$.
Indeed, with respect to the fundamental question how far even despite the presence of such saturation effects,
concentration phenomena in the sense of spontaneous emergence of
locally unbounded densities may occur in (\ref{PE}), when considered under homogeneous Neumann boundary conditions 
in smoothly bounded $N$-dimensional domains with $N\ge 2$,
the number $\alpha_c(N):=\frac{N-2}{N-1}$ has been found to play the role of a critical exponent:
Whenever $F$ is suitably smooth and such that
\be{F0}
	|F(\xi)| \le \Kf \cdot (\xi+1)^{-\frac{\alpha}{2}}
	\qquad \mbox{for all } \xi\ge 0
\ee
with some $\Kf>0$ and $\alpha>\alpha_c$, an associated initial-boundary value problem admits globally bounded classical solutions
for widely arbitrary initial data, while if $N\ge 3$, $\Omega$ is a ball and
\be{fsuper}
	F(\xi) \ge \kf \cdot (\xi+1)^{-\frac{\alpha}{2}}
	\qquad \mbox{for all } \xi\ge 0
\ee
with some $\kf>0$ and $\alpha<\alpha_c$, then some radially symmetric classical solutions exist which blow up at some finite time
$T>0$ in the sense that $\limsup_{t\nearrow T} \|n(\cdot,t)\|_{L^\infty(\Omega)} =\infty$
(\cite{win_ct_flim_critexp}; cf.~also \cite{negreanu_tello} and \cite{win_ct_unify}).\abs
Having this fairly complete knowledge for (\ref{PE}) at hand, in this study we intend to accomplish an apparently natural next step
toward an understanding of models for tactic movement in realistic situations:
Namely, we plan to examine to which extent the above blow-up related dichotomy may be affected if buoyancy-induced interaction
of a considered population with its liquid environment is additionally accounted for.
In fact, experimental evidence indicates a substantial relevance of chemotaxis-fluid couplings in several contexts
(\cite{coll}, \cite{deshmane}, \cite{miller}, \cite{taub}, \cite{goldstein}),
and both numerical (\cite{goldstein},  \cite{lorz_CMS}) and some rigorous analytical precedents 
(\cite{kiselev_ryzhik1}, \cite{kiselev_ryzhik2}, \cite{kiselev_xu}, \cite{he_tadmor}, \cite{evje_win}) 
witness significant quantitative effects
of fluid interaction in various specific situations, mostly involving a given and hence passive fluid.\abs
To take aim at this question in a fairly general three-dimensional setting compatible with the seminal modeling approach
in \cite{goldstein}, we shall below consider a fully parabolic counterpart of (\ref{PE}), coupled to the 
incompressible Stokes equations, in the framework of the initial-boundary value problem
\be{0}
	\left\{ \begin{array}{lcll}
	n_t + u\cdot\na n &=& \Delta n - \na \cdot \Big( n F(|\na c|^2) \na c\Big),
	\qquad & x\in \Omega, \ t>0, \\[1mm]
	c_t + u\cdot\na c &=& \Delta c - c + n,
	\qquad & x\in \Omega, \ t>0, \\[1mm]
	u_t &=& \Delta u + \na P + n\na\Phi,
	\qquad \na \cdot u=0,
	\qquad & x\in \Omega, \ t>0, \\[1mm]
	& & \hs{-29mm}
	\frac{\partial n}{\partial\nu}=\frac{\partial c}{\partial\nu}=0, \quad u=0,
	\qquad & x\in \pO, \ t>0, \\[1mm]
	& & \hs{-29mm}
	n(x,0)=n_0(x), \quad c(x,0)=c_0(x), \quad u(x,0)=u_0(x),
	\qquad & x\in\Omega,
	\end{array} \right.
\ee
where $\Omega\subset\R^3$ is a bounded domain with smooth boundary, where $\Phi\in W^{2,\infty}(\Om)$ represents a given
gravitational potential, and where $F$ appropriately generalizes the choice in (\ref{fp}).
Since any exploding solution of the corresponding fluid-free chemotaxis-only system trivially extends to an accordingly singular
solution of (\ref{0}) upon letting $\Phi\equiv const.$ and $u_0\equiv 0$ therein, 
and since in view of the above criticality feature of $\alpha_c(3)=\frac{1}{2}$
we thus expect blow-up to occur in (\ref{0}) whenever
$F$ satisfies (\ref{fsuper}) with some $\kf>0$ and $\alpha < \frac{1}{2}$,
the major part of our ambition seems to consist in deciding whether or not the assumption (\ref{F0})
with $\Kf>0$ and arbitrary $\alpha>\frac{1}{2}$ remains sufficient to ensure global smooth solvability and 
boundedness also in the context of the full problem (\ref{0}).\abs
{\bf Conditional $L^\infty$ estimates for signal gradients in chemotaxis-Stokes systems. Main results I.} 
In order to develop an approach which is not only appropriate to address this topic, but which 
is suitably broad so as to bear some potential for meaningful applications in different frameworks beyond the particular present one,
prior to our analysis of the specific system (\ref{0}) we shall take up the general problem of efficiently controlling
signal gradients in contexts of chemotaxis processes coupled to Stokes fluids.
This problem has formed a natural core of several precedent studies concerned with models related to (\ref{0}),
and a considerable variety of respective ad-hoc approaches, based either on the detection of certain entropy-like features, 
or on maximal Sobolev regularity results, 
or on heat semigroup estimates,
can be found in the literature
(cf.~\cite{liu_wangyifu},
\cite{wangyulan_xiang_JDE2015}, \cite{wangyulan_xiang_JDE2016}, \cite{cao_ZAMP}, \cite{win_CVPDE} for examples).\abs
An essential characteristic of our subsequent analysis will now be made up by the intention to propose
the use of certain conditional $L^\infty$ estimates for taxis gradients as a possible alternative to previous strategies.
The derivation of such bounds in suitably general settings will form a fairly extensive first part of this manuscript,
but as a considerable benefit of the related efforts, inter alia indicating their potential value beyond the context of (\ref{0}),
we will thereafter be able to conveniently accomplish our original objective by means of quite a compressed argument.\abs
In order to make this more precise, 
in a bounded domain $\Omega\subset\R^3$ and with given $T\in (0,\infty]$, $m>0, K>0$ and $\vartheta\in (\frac{3}{4},1)$, 
let us consider solutions $(z,v,Q)$, with
\be{vz}
	\left\{ \begin{array}{l}
	z\in \bigcap_{q>3} C^0([0,T);W^{1,q}(\Om))	 \cap C^{2,1}(\bom\times (0,T)), \\[1mm]
	v\in 	
	C^0([0,T);D(A^\vartheta)) \cap C^{2,1}(\bom\times (0,T);\R^3) \qquad \mbox{and} \\[1mm]
	Q \in C^{1,0}(\Om\times (0,T)),
	\end{array} \right.
\ee
of the problem
\be{0vz}
	\left\{ \begin{array}{lcll}
	z_t + v\cdot\na z &=& \Delta z - z + f,
	\qquad & x\in \Omega, \ t\in (0,T), \\[1mm]
	v_t &=& \Delta v + \na Q + f\na \Phi,
	\qquad \na \cdot v=0,
	\qquad & x\in \Omega, \ t\in (0,T), \\[1mm]
	& & \hs{-28mm}
	\frac{\partial z}{\partial\nu}=0, \quad v=0,
	\qquad & x\in \pO, \ t\in (0,T), \\[1mm]
	& & \hs{-28mm}
	z(x,0)=z_0(x), \quad v(x,0)=v_0(x),
	\qquad & x\in\Omega,
	\end{array} \right.
\ee
where $z_0, v_0$ and $f$ are given functions which are such that
\be{init}
	\left\{ \begin{array}{l}
	z_0\in W^{1,\infty}(\Om) \mbox{ is nonnegative\quad  and \quad }
	v_0\in D(A^\vartheta) \qquad \mbox{with} \\[1mm]
	\|z_0\|_{W^{1,\infty}(\Om)} \le K \quad \mbox{and} \quad \|A^\vartheta v_0\|_{L^2(\Om)} \le K,
	\end{array} \right.
\ee
and that
\be{f}
	\left\{ \begin{array}{l}
	f\in C^0(\bom\times [0,T))\mbox{ is nonnegative with } \\[1mm]
	\io f(\cdot,t) \le m
	\quad \mbox{for all } t\in (0,T).
	\end{array} \right.
\ee
Here, as throughout the sequel, 
we let $A=-\proj \Delta$ and $(A^\vartheta)_{\vartheta>0}$ denote the realization of the 
Stokes operator in $L^2(\Omega;\R^3)$, with domain given by
$D(A):=W^{2,2}(\Om;\R^3)\cap W_{0,\sigma}^{1,2}(\Om)$, and the family of its corresponding fractional powers,
respectively, where $\proj$ represents the Helmholtz projection on $L^2(\Om;\R^3)$,
and where 
$W_{0,\sigma}^{1,2}(\Omega):=W_0^{1,2}(\Omega;\R^3)\cap L^2_\sigma(\Om)$
with $L^2_\sigma(\Omega):=\{\varphi\in L^2(\Omega;\R^3) \ | \ \na\cdot\varphi=0\}$.\abs
Within this framework, by means of quite a direct application of standard regularity estimates for the Stokes semigroup
we shall firstly obtain the following statement on regularity of the corresponding fluid field,
conditional in that some temporally independent $L^p$ bounds for the forcing term $f$ in (\ref{0vz}) are involved.
\begin{prop}\label{prop1}
  Let $\Omega\subset\R^3$ be a bounded domain with smooth boundary, and
  suppose that $\vartheta\in (\frac{3}{4},1)$, $m>0, K>0$, $p\ge 2$, $\theta\in (\frac{1}{4},\vartheta]$ and $\eta>0$.
  Then there exists $C(m,K,p,\theta,\eta)>0$ such that any $f,z,v$ and $Q$ fulfilling (\ref{f}), (\ref{vz}) and (\ref{0vz}),
  with some $T\in (0,\infty]$ and some $(z_0,v_0)$ complying with (\ref{init}), satisfies
  \be{1.1}
	\|A^\theta v(\cdot,t)\|_{L^2(\Om)} 
	\le C(m,K,p,\theta,\eta) \cdot \bigg\{ 1 + \sup_{s\in (0,t)} \|f(\cdot,s)\|_{L^p(\Om)} \bigg\}^
		{\frac{p}{p-1} \cdot (\frac{4\theta-1}{6} + \eta)} 
	\qquad \mbox{for all } t\in (0,T).
  \ee
\end{prop}
On the basis thereof and of appropriate regularization features of Neumann heat semigroups, we shall secondly derive
uniform estimates for the first solution components in (\ref{0vz}) and their gradients, subject to a conditionality
property similar to that expressed in (\ref{1.1}):
\begin{theo}\label{theo2}
  Suppose that $\Omega\subset\R^3$ is a bounded domain with smooth boundary, and
  let $\vartheta\in (\frac{3}{4},1)$, $m>0, K>0$, $p>3$ and $\eta>0$.
  Then one can find $C(m,K,p,\eta)>0$ with the property that given any $T\in (0,\infty]$ as well as functions $f,z,v$ and $Q$ 
  which satisfy (\ref{f}), (\ref{vz}) and (\ref{0vz}) with some $(z_0,v_0)$ fulfilling (\ref{init}), we have
  \be{2.1}
	\|z(\cdot,t)\|_{W^{1,\infty}(\Om)}
	\le C(m,K,p,\eta) \cdot \bigg\{ 
	1+ \sup_{s\in (0,t)} \|f(\cdot,s)\|_{L^p(\Om)} \bigg\}^{\frac{p}{p-1} \cdot (\frac{2}{3}+\eta)} 
	\qquad \mbox{for all } t\in (0,T).
  \ee
\end{theo}
{\bf Approaching criticality in boundedness results for (\ref{0}). Main results II.} \quad
With the above preparations at hand, in Section \ref{sect3} we can thereafter return to our original question.
In fact, Proposition \ref{prop1} and, especially, Theorem \ref{theo2} will enable us to make sure
by means of a fairly condensed argument
that global existence and boundedness in (\ref{0}) can indeed be achieved under the assumption (\ref{F0})
within the entire range of $\alpha$ known as the regime essentially maximal in this regard for (\ref{0}):
\begin{theo}\label{theo11}
  Let $\Omega\subset\R^3$ be a bounded domain with smooth boundary, and let $\Phi\in W^{2,\infty}(\Om)$ and $F\in C^2([0,\infty))$
  be such that
  \be{F}
	|F(\xi)| \le \Kf \cdot (1+\xi)^{-\frac{\alpha}{2}}
	\qquad \mbox{for all } \xi\ge 0
  \ee
  with some $\Kf>0$ and
  \bas
	\alpha>\frac{1}{2}.
  \eas
  Then for any
  \be{incu}
	\left\{ \begin{array}{l}
	n_0\in C^0(\bom) \mbox{ such that $n_0\ge 0$}, \\[1mm]
	c_0 \in W^{1,\infty}(\Om) \mbox{ such that } c_0\ge 0 \qquad \mbox{and} \\[1mm]
	u_0 \in \bigcup_{\vartheta\in (\frac{3}{4},1)} D(A^\vartheta),
	\end{array} \right.
  \ee
  there exist uniquely determined
  \bas
	\left\{ \begin{array}{l}
	n\in C^0(\bom\times [0,\infty)) \cap C^{2,1}(\bom\times (0,\infty)), \\[1mm]
	c\in \bigcap_{q>3} C^0([0,\infty);W^{1,q}(\Om)) \cap C^{2,1}(\bom\times (0,\infty)) \qquad \mbox{and} \\[1mm]
	u\in \bigcup_{\vartheta\in (\frac{3}{4},1)} C^0([0,\infty);D(A^\vartheta)) \cap C^{2,1}(\bom\times (0,\infty);\R^3)
	\end{array} \right.
  \eas
  such that $n\ge 0$ and $c\ge 0$ in $\Omega\times (0,\infty)$, and that with some $P\in C^{1,0}(\Om\times (0,\infty))$, the collection
  $(n,c,u,P)$ solves (\ref{0}) in the classical sense in $\Om\times (0,\infty)$.
  Moreover, this solution is bounded in the sense that with some $\vartheta\in (\frac{3}{4},1)$ and $C>0$,
  \be{11.1}
	\|n(\cdot,t)\|_{L^\infty(\Om)}
	+ \|c(\cdot,t)\|_{W^{1,\infty}(\Om)}
	+ \|A^\vartheta u(\cdot,t)\|_{L^2(\Om)}
	\le C
	\qquad \mbox{for all } t>0.
  \ee
\end{theo}

\vs{5mm}
Before going into details, let us once again emphasize that
applications of the above general estimates, and especially that from Theorem \ref{theo2}, to much a larger variety of 
chemotaxis-Stokes systems beyond the particular context of (\ref{0}) are conceivable.
Although in order to keep this manuscript reasonably focused we do not discuss this here in detail, we at least may
announce that one further example concerned with a quasilinear Keller-Segel system coupled to the Stokes equations 
will be addressed in \cite{win_AML}.
\mysection{Conditional bounds. Proofs of Proposition \ref{prop1} and Theorem \ref{theo2}}\label{sect2}
\subsection{Fluid regularity}
Let us first employ known smoothing properties of the Stokes semigroup to see that in the considered context of 
solutions $(z,v,Q)$ of (\ref{0vz}),
regularity properties of the respective fluid can be related to the respective forcing term through
the expressions
\be{M}
	M_p(t):=1 + \sup_{s\in (0,t)} \|f(\cdot,s)\|_{L^p(\Om)},
	\qquad t\in (0,T),
\ee
well-defined due to (\ref{f}) actually for all $p\ge 1$, 
in the quantitative manner specified in Proposition \ref{prop1}:\abs
\proofc of Proposition \ref{prop1}. \quad
  Since $\theta<1$, we may assume that $\eta<\frac{2(1-\theta)}{3}$, which ensures that $7-4\theta-6\eta>3$, whence observing that
  $7-4\theta<6$ due to the restriction $\theta>\frac{1}{4}$, we obtain that
  \bas
	1<\frac{6}{7-4\theta} < \frac{6}{7-4\theta-6\eta} <2.
  \eas
  As, on the other hand, we are presupposing that $p\ge 2$, if we choose any
  \be{1.3}
	\lam=\lam(p,\theta,\eta) \in \Big(\frac{6}{7-4\theta} \, , \, \frac{6}{7-4\theta-6\eta} \Big],
  \ee
  then necessarily $1<\lam < 2\le p$, so that by continuity of the Helmholtz projection on $L^\lam(\Om;\R^3)$
  (\cite{fujiwara_morimoto}), we can find $C_1=C_1(p,\theta,\eta)>0$ such that abbreviating
  $C_2:=\|\na\Phi\|_{L^\infty(\Om)}$ and $a:=\frac{p(\lam-1)}{(p-1)\lam}$, 
  due to the H\"older inequality we have
  \bea{1.32}
	\big\| \proj [\varphi\na\Phi]\big\|_{L^\lam(\Om)}
	&\le& C_1 \|\varphi\na\Phi\|_{L^\lam(\Om)} \nn\\
	&\le& C_1 C_2 \|\varphi\|_{L^\lam(\Om)} \nn\\
	&\le& C_1 C_2 \|\varphi\|_{L^p(\Om)}^a \|\varphi\|_{L^1(\Om)}^{1-a} 
	\qquad \mbox{for all } \varphi\in C^0(\bom).
  \eea
  We next use that $\theta\le \vartheta$ and
  invoke known regularization features of the Stokes semigroup $(e^{-tA})_{t\ge 0}$ (\cite[p.201]{giga1986}) to infer that
  with some $C_3=C_3(\theta)>0$, $C_4=C_4(p,\theta,\eta)>0$ and $\mu>0$,
  \be{1.33}
	\|A^\theta e^{-tA} \varphi\|_{L^2(\Om)}
	\le C_3 \|A^\vartheta \varphi\|_{L^2(\Om)}
	\qquad \mbox{for all $t\ge 0$ and } \varphi\in D(A^\vartheta)
  \ee
  as well as
  \be{1.34}
	\|A^\theta e^{-tA} \varphi\|_{L^2(\Om)}
	\le C_4 t^{-\theta-\frac{3}{2}(\frac{1}{\lam}-\frac{1}{2})} e^{-\mu t} \|\varphi\|_{L^\lam(\Om)}
	\qquad \mbox{for all $t> 0$ and } \varphi\in C^0(\bom;\R^3).
  \ee
  We now let $m>0$ and $K>0$ be given, and suppose that $f,z_0, v_0, v,z$ and $Q$ satisfy 
  (\ref{vz}), (\ref{0vz}), (\ref{f}) and (\ref{init}) with some $T\in (0,\infty]$.
  Then taking $(M_p(t))_{t\in (0,T)}$ as accordingly introduced in (\ref{M}), from (\ref{1.32}) we obtain that
  \bas
	\big\| \proj [f(\cdot,s)\na\Phi]\big\|_{L^\lam(\Om)}
	\le C_1 C_2 m^{1-a} M_p^a(t)
	\qquad \mbox{for all $t\in (0,T)$ and any } s\in (0,t),
  \eas
  whence a combination of (\ref{1.33}) with (\ref{1.34}) shows that
  \bea{1.4}
	\|A^\theta v(\cdot,t)\|_{L^2(\Om)}
	&=& \bigg\| A^\theta e^{-tA} v_0
	+ \int_0^t A^\theta e^{-(t-s)A} \proj [ f(\cdot,s)\na \Phi] ds \bigg\|_{L^2(\Om)} \nn\\
	&\le& C_3 \|A^\vartheta v_0\|_{L^2(\Om)}
	+ C_4 \int_0^t (t-s)^{-\theta-\frac{3}{2}(\frac{1}{\lam}-\frac{1}{2})} e^{-\mu (t-s)} 
		\big\|\proj [f(\cdot,s)\na\Phi] \big\|_{L^\lam(\Om)} ds \nn\\
	&\le& C_3 K
	+ C_1 C_2 C_4 C_5 m^{1-a} M_p^a(t)
	\qquad \mbox{for all } t\in (0,T)
  \eea  
  thanks to (\ref{init}), where 
  $C_5=C_5(p,\theta,\eta):=\int_0^\infty \sigma^{-\theta-\frac{3}{2}(\frac{1}{\lam}-\frac{1}{2})} e^{-\mu\sigma} d\sigma$
  is finite due to the lower bound for $\lam$ implied by (\ref{1.3}), which namely warrants that
  \bas
	\theta+\frac{3}{2}\cdot\Big(\frac{1}{\lam}-\frac{1}{2}\Big)
	< \theta + \frac{3}{2} \cdot \Big(\frac{7-4\theta}{6}-\frac{1}{2}\Big) =1.
  \eas
  Since, conversely, the upper bound for $\lam$ in (\ref{1.3}) guarantees that
  \bas
	a=\frac{p}{p-1} \cdot \Big(1-\frac{1}{\lam}\Big)
	\le \frac{p}{p-1} \cdot \Big(1-\frac{7-4\theta-6\eta}{6}\Big) 
	= \frac{p}{p-1} \cdot \Big(\frac{4\theta-1}{6}+\eta\Big),
  \eas
  and that thus the inequality $M_p\ge 1$ ensures that
  $M_p^a(t) \le M_p^{\frac{p}{p-1}\cdot (\frac{4\theta-1}{6}+\eta)}(t)$ for all $t\in (0,T)$,
  we conclude that (\ref{1.1}) is a consequence of (\ref{1.4}).
\qed
In order to appropriately prepare an application of the latter in the course of our derivation of Theorem \ref{theo2},
we recall known embedding properties enjoyed by the domains of the operators $A^\theta$ appearing in
(\ref{1.1}) to turn the latter into integral bounds for the fluid velocity itself:
\begin{cor}\label{cor4}
  Let $\vartheta\in (\frac{3}{4},1)$, and let $m>0, K>0, p\ge 2, r>3$ and $\eta>0$.
  Then there exists $C(m,K,p,r,\eta)>0$ such that if (\ref{f}) and (\ref{init}) as well as (\ref{vz}) and (\ref{0vz}) hold with some
  $T\in (0,\infty]$, it follows that
  \be{4.1}
	\|v(\cdot,t)\|_{L^r(\Om)} \le
	C(m,K,p,r,\eta) \cdot \bigg\{ 1 + \sup_{s\in (0,t)} \|f(\cdot,s)\|_{L^p(\Om)} \bigg\}^
		{\frac{p}{p-1} \cdot (\frac{r-3}{3r}+\eta)}
	\qquad \mbox{for all } t\in (0,T).
  \ee
\end{cor}
\proof
  Since $4\vartheta-3\ge 0$ and thus $\frac{(4\vartheta-3)r+6}{5r}>0$, we may restrict ourselves to considering the case when
  $\eta \le \frac{(4\vartheta-3)r+6}{5r}$, in which
  \be{4.2}
	\theta=\theta(r,\eta):=\frac{3r-6+5r\eta}{4r}
  \ee
  satisfies
  \bas
	\theta \le \frac{3r-6+5r \cdot \frac{(4\vartheta-3)r +6}{5r}}{4r} =\vartheta.
  \eas
  As (\ref{4.2}) furthermore clearly ensures that
  \be{4.3}
	\theta>\frac{3r-6}{4r},
  \ee
  and that thus
  \bas
	\theta>\frac{3-\frac{6}{r}}{4}>\frac{1}{4}
  \eas
  due to our assumption that $r>3$, 
  we may employ Proposition \ref{prop1} to find $C_1=C_1(m,K,p,r,\eta)$ with the property that
  if $T\in (0,\infty]$ as well as (\ref{f}), (\ref{init}), (\ref{vz}) and (\ref{0vz}) hold, then
  with $(M_p(t))_{t\in (0,T)}$ taken from (\ref{M}) we have
  \be{4.33}
	\|A^\theta v(\cdot,t)\|_{L^2(\Om)}
	\le C_1 M_p^{\frac{p}{p-1} \cdot \frac{4\theta-1+\eta}{6}}(t)
	\qquad \mbox{for all } t\in (0,T).
  \ee
  Apart from that, (\ref{4.3}) implies that
  \bas
	2\theta-\frac{3}{2}+\frac{3}{r} > 2\cdot \frac{3r-6}{4r} - \frac{3}{2} + \frac{3}{r}=0,
  \eas
  from which 
  it follows that $D(A^\theta)$ is continuously embedded into $L^r(\Omega;\R^3)$ (\cite{giga1981_the_other}, \cite{henry}), and that hence we can
  pick $C_2=C_2(r,\eta)>0$ such that
  \be{4.4}
	\|\varphi\|_{L^r(\Om)} \le C_2 \|A^\theta \varphi\|_{L^2(\Om)}
	\qquad \mbox{for all } \varphi\in D(A^\theta).
  \ee
  It now only remains to observe that if we 
  assume (\ref{f}), (\ref{init}), (\ref{vz}) and (\ref{0vz})
  to be satisfied with some $T\in (0,\infty]$, then a combination of (\ref{4.4}) with (\ref{4.33}) 
  shows that with $M_p$ as accordingly defined through (\ref{M}) we have
  \bas
	\|v(\cdot,t)\|_{L^r(\Om)}
	\le C_1 C_2 M_p^{\frac{p}{p-1} \cdot \frac{4\theta-1+\eta}{6}}(t)
	\qquad \mbox{for all } t\in (0,T),
  \eas
  and that 
  \bas
	\frac{4\theta-1+\eta}{6} 
	= \frac{4\cdot \frac{3r-6+5r\eta}{4r}-1+\eta}{6}
	= \frac{2r-6+6r\eta}{6r}
	= \frac{r-3}{3r} + \eta
  \eas
  by (\ref{4.2}).
\qed
\subsection{Uniform estimates for taxis gradients}
The first step toward our analysis concerned with the first equation in (\ref{0vz}) acts at levels which,
through relying on a standard zero-order testing procedure only, 
do not explicitly rely on any propertiy of the considered fluid field beyond its mere solenoidality.
Although the estimates for $z$ thereby obtained involve topologies which are yet quite far from those to be finally addressed,
they will form an essential basic information in the course of an interpolation argument performed in Lemma \ref{lem6} below.
\begin{lem}\label{lem5}
  Let $\vartheta\in (\frac{3}{4},1)$. Then for any $m>0, K>0, p\ge 2$ and $q\ge 2$ there exists $C(m,K,p,q)>0$ such that if
  (\ref{f}), (\ref{init}), (\ref{vz}) and (\ref{0vz}) hold with some $T\in (0,\infty]$, then
  \be{5.1}
	\|z(\cdot,t)\|_{L^q(\Om)}
	\le C(m,K,p,q) \cdot \bigg\{ 1 + \sup_{s\in (0,t)} \|f(\cdot,s)\|_{L^p(\Om)} \bigg\}^{\frac{p}{p-1} \cdot \frac{q-1}{3q}}
	\qquad \mbox{for all } t\in (0,T).
  \ee
\end{lem}
\proof
  By continuity of the embedding $W^{1,2}(\Om) \hra L^6(\Om)$, we can fix $C_1>0$ such that
  \be{5.2}
	\|\varphi\|_{L^6(\Om)}^2 \le C_1 \io |\na\varphi|^2 + C_1 \io \varphi^2
	\qquad \mbox{for all } \varphi\in W^{1,2}(\Om).
  \ee
  Therefore, assuming (\ref{f}), (\ref{init}), (\ref{vz}) and (\ref{0vz}) with some $T\in (0,\infty]$, in the identity
  \be{5.3}
	\frac{1}{q} \frac{d}{dt} \io z^q
	+ \frac{4(q-1)}{q^2} \io |\na z^\frac{q}{2}|^2 + \io z^q
	= \io z^{q-1} f,
	\qquad t\in (0,T),
  \ee
  as obtained upon testing the first equation in (\ref{0vz}) by $z^{q-1}$ due to the solenoidality of $v$, 
  we may combine the H\"older inequality with (\ref{5.2})
  and Young's inequality to estimate
  \bea{5.4}
	\io z^{q-1} f
	&=& \io (z^\frac{q}{2})^\frac{2(q-1)}{q} f \nn\\
	&\le& \|z^\frac{q}{2}\|_{L^6(\Om)}^\frac{2(q-1)}{q} \|f\|_{L^\frac{3q}{2q+1}(\Om)} \nn\\
	&\le& C_1^\frac{q-1}{q} \cdot \bigg\{ \io |\na z^\frac{q}{2}|^2 + \io z^q \bigg\} \cdot \|f\|_{L^\frac{3q}{2q+1}(\Om)} \nn\\
	&\le& \frac{q-1}{q} \cdot \bigg\{ \io |\na z^\frac{q}{2}|^2 + \io z^q \bigg\}
	+ \frac{1}{q} \cdot \Big\{ C_1^\frac{q-1}{q} \|f\|_{L^\frac{3q}{2q+1}(\Om)} \Big\}^q \nn\\
	&=& \frac{q-1}{q} \io |\na z^\frac{q}{2}|^2 
	+ \frac{q-1}{q} \io z^q
	+ \frac{C_1^{q-1}}{q} \|f\|_{L^\frac{3q}{2q+1}(\Om)}^q
	\qquad \mbox{for all } t\in (0,T).
  \eea
  Here since $1<\frac{3q}{2q+1}<\frac{3}{2}<p$, we may once again rely on the H\"older inequality to see that thanks to (\ref{f}),
  with $(M_p(t))_{t\in (0,T)}$ as in (\ref{M}) we have
  \bas
	\|f\|_{L^\frac{3q}{2q+1}(\Om)}^q
	&\le& \|f\|_{L^p(\Om)}^\frac{p(q-1)}{3(p-1)} \|f\|_{L^1(\Om)}^\frac{2pq+p-3q}{3(p-1)} \\
	&\le& m^\frac{2pq+p-3q}{3(p-1)} M_p^\frac{p(q-1)}{3(p-1)}(t)
	\qquad \mbox{for all } t\in (0,T).
  \eas
  As
  \bas
	\frac{4(q-1)}{q^2} \ge \frac{q-1}{q}
	\qquad \mbox{and} \qquad
	1-\frac{q-1}{q} \ge \frac{1}{2}\ge \frac{1}{q}
  \eas
  due to the fact that $q\ge 2$, from (\ref{5.3}) and (\ref{5.4}) we therefore obtain that for any choice of $t_0\in (0,T)$,
  \bas
	\frac{d}{dt} \io z^q + \io z^q
	\le C_1^{q-1} m^\frac{2pq+p-3q}{3(p-1)} M_p^\frac{p(q-1)}{3(p-1)}(t_0)
	\qquad \mbox{for all } t\in (0,t_0),
  \eas
  and that thus, according to an ODE comparison argument,
  \bas
	\io z^q
	\le \max \bigg\{ \io z_0^q \, , \, C_1^{q-1} m^\frac{2pq+p-3q}{3(p-1)} M_p^\frac{p(q-1)}{3(p-1)}(t_0) \bigg\}
	\qquad \mbox{for all } t\in (0,t_0].
  \eas
  When evaluated at $t=t_0$, in view of (\ref{init}) this readily yields the claim.
\qed
Indeed, the latter can be used to derive the following preliminaryl information 
on control of spatial $W^{1,\infty}$ norms of $z$, as addressed in Theorem \ref{theo2}.
In its formulation and throughout the remainder of this section, we let $B=B_q$ denote the sectorial realization of $-\Delta+1$
under homogeneous Neumann boundary conditions in $\bigcap_{q>1} L^q(\Omega)$, and let $(B^\beta)_{\beta>0}$
represent 
the associated family of positive fractional powers $B^\beta=B_q^\beta$.
Then the respective domains $D(B_q^\beta)$ are continuously embedded into $W^{1,\infty}(\Om)$ whenever $2\beta-\frac{3}{q}>1$,
whence for such parameters the quantities under consideration can be estimated against expressions of the form
$\|B^\beta z\|_{L^q(\Om)}$, together with the Lebesgue norms from Lemma \ref{lem5}, through interpolation. 
This is substantiated in the following statement which circumvents problems potentially resulting from possibly insufficient
regularity properties of $z_0$, as merely assumed here to belong to $W^{1,\infty}(\Omega)$ but not necessarily to $D(B_q^\beta)$
for any $\beta>\frac{1}{2}$, by subtracting a suitable correction.
\begin{lem}\label{lem6}
  Let $\vartheta\in (\frac{3}{4},1)$, and let $m>0, K>0, p\ge 2, q\ge 2$ and $\beta\in (\frac{1}{2},1)$ be such that
  \be{6.1}
	q\cdot (2\beta-1)>3.
  \ee
  Then for all $\eta>0$ one can find $C(m,K,p,q,\beta,\eta)>0$ such that whenever
  (\ref{f}), (\ref{init}), (\ref{vz}) and (\ref{0vz}) are satisfied with some $T\in (0,\infty]$, we have
  \bea{6.2}
	\|z(\cdot,t)-e^{-tB} z_0\|_{W^{1,\infty}(\Om)} 
	&\le& C(m,K,p,q,\beta,\eta) \cdot 
	\bigg\{ 1 + \sup_{s\in (0,t)} \|f(\cdot,s)\|_{L^p(\Om)} \bigg\}^{\frac{p}{p-1} \cdot \frac{(q-1)(2q\beta-q-3)}{6q^2\beta}}
	\times \nn\\
	& & \hs{29mm} \times \,
	\bigg\{ 1 + \Big\| B^\beta \Big(z(\cdot,t)-e^{-tB} z_0\Big) \Big\|_{L^q(\Om)} \bigg\}^{\frac{q+3}{2q\beta}+\eta)}
  \eea
  for all $t\in (0,T)$.
\end{lem}
\proof
  As (\ref{6.1}) warrants that
  \bas
	\frac{q+3}{3q} = \frac{1}{2} + \frac{3}{2q} < \frac{1}{2} + \frac{3}{2} \cdot \frac{2\beta-1}{3}=\beta,
  \eas
  we may assume that $\eta$ is so small that
  \be{6.3}
	\gamma=\gamma(q,\eta):=\frac{q+3}{2q} + \beta\eta
  \ee
  satisfies $\gamma<\beta$.
  We may therefore draw on a standard interpolation property enjoyed by fractional powers of sectorial operators in quite
  general settings (see e.g.~\cite[Theorem 2.14.1]{friedman}) in choosing $C_1=C_1(q,\beta,\eta)>0$ such that
  \be{6.4}
	\|B^\gamma \varphi\|_{L^q(\Om)} 
	\le C_1 \|B^\beta\varphi\|_{L^q(\Om)}^\frac{\gamma}{\beta} \|\varphi\|_{L^q(\Om)}^\frac{\beta-\gamma}{\beta}
	\qquad \mbox{for all } \varphi\in D(B^\beta).
  \ee
  Furthermore, the evient fact that $\gamma>\frac{q+3}{2q}$ guarantees that $2\gamma<\frac{3}{q}>1$, meaning that
  $D(B^\gamma)$ is continuously embedded into $W^{1,\infty}(\Om)$ (\cite{henry}), and that thus there exists
  $C_2=C_2(q,\beta,\eta)>0$ fulfilling
  \be{6.5}
	\|\varphi\|_{W^{1,\infty}(\Om)} \le C_2 \|B^\gamma \varphi\|_{L^q(\Om)}
	\qquad \mbox{for all } \varphi\in D(B^\gamma).
  \ee
  As a final preparation, we take $C_3=C_3(q)>0$ in such a way that
  \be{6.55}
	\|\varphi\|_{L^q(\Om)} \le C_3 \|\varphi\|_{W^{1,\infty}(\Om)}
	\qquad \mbox{for all } \varphi\in W^{1,\infty}(\Om),
  \ee
  and invoke Lemma \ref{lem5} to see that whenever $m>0, K>0$ and $p\ge 2$, we can fix $C_4=C_4(m,K,p,q)>0$
  with the property that if (\ref{f}), (\ref{init}), (\ref{vz}) and (\ref{0vz}) hold with some $T\in (0,\infty]$, then
  \be{6.6}
	\|z(\cdot,t)\|_{L^q(\Om)}
	\le C_4 M_p^{\frac{p}{p-1} \cdot \frac{q-1}{3q}}(t)
	\qquad \mbox{for all } t\in (0,T),
  \ee
  again with $(M_p(t))_{t\in (0,T)}$ as defined in (\ref{M}).\abs
  Thus, given any such $T$ and arbitrary $f,z_0,v_0,z,v$ and $Q$ satisfying (\ref{f}), (\ref{init}), (\ref{vz}) and (\ref{0vz}),
  abbreviating 
  \be{6.66}
	\wh{z}(\cdot,t):=z(\cdot,t)-e^{-tB} z_0,
	\qquad t\in (0,T),
  \ee
  we can firstly rely on (\ref{6.6}), (\ref{6.55}) and (\ref{init}) to see that 
  with the correspondingly defined function $M_p \ge 1$, since $e^{-tB}$ is nonexpansive on $L^q(\Om)$ for
  all $t>0$ we have
  \bas
	\|\wh{z}(\cdot,t)\|_{L^q(\Om)}
	&\le& \|z(\cdot,t)\|_{L^q(\Om)}
	+ \|e^{-tB} z_0\|_{L^q(\Om)} \\
	&\le& C_4 M_p^{\frac{p}{p-1} \cdot \frac{q-1}{q}}(t) 
	+ \|z_0\|_{L^q(\Om)} \\
	&\le& C_4 M_p^{\frac{p}{p-1} \cdot \frac{q-1}{q}}(t) 
	+ C_3 K \\
	&\le& C_5 M_p^{\frac{p}{p-1} \cdot \frac{q-1}{3q}}(t)
	\qquad \mbox{for all } t\in (0,T),
  \eas
  where $C_5=C_5(m,K,p,q):=\max \{ C_4 \, , \, C_3 K\}$.
  Combining (\ref{6.5}) with (\ref{6.4}) hence shows that
  \bas
	\|\wh{z}(\cdot,t)\|_{W^{1,\infty}(\Om)}
	&\le& C_2 \|B^\gamma \wh{z}(\cdot,t)\|_{L^q(\Om)} \\
	&\le& C_1 C_2 \|B^\beta \wh{z}(\cdot,t)\|_{L^q(\Om)}^\frac{\gamma}{\beta} 
		\|\wh{z}(\cdot,t)\|_{L^q(\Om)}^\frac{\beta-\gamma}{\beta} \\
	&\le& C_1 C_2 C_5^\frac{\beta-\gamma}{\beta} M_p^{\frac{p}{p-1} \cdot \frac{q-1}{3q} \cdot \frac{\beta-\gamma}{\beta}}(t)
		\|B^\beta \wh{z}(\cdot,t)\|_{L^q(\Om)}^\frac{\gamma}{\beta}
	\qquad \mbox{for all } t\in (0,T),
  \eas
  from which (\ref{6.2}) results upon observing that
  \bas
	\frac{\gamma}{\beta}=\frac{q-3}{2q\beta} + \eta
  \eas
  by (\ref{6.3}), and that
  \bas
	\frac{p}{p-1} \cdot \frac{q-1}{3q} \cdot \frac{\beta-\gamma}{\beta}
	&\le& \frac{p}{p-1} \cdot \frac{q-1}{3q} \cdot \frac{\beta-\frac{q+3}{2q}}{\beta} \\
	&=& \frac{p}{p-1} \cdot \frac{(q-1)(2q\beta-q+3)}{6q^2\beta}
  \eas
  due to the inequality $\gamma\ge \frac{q+3}{2q}$.
\qed
The core of our analysis in this section is now contained in the following estimate for the higher-order expressions on the right
of (\ref{6.2}) against the quantities in (\ref{M}).
This is achieved by means of smoothing estimates for the semigroup $(e^{-tB})_{t\ge 0}$, where the corresponding
convective contribution to the first equation in (\ref{0vz}) will be treated as a lower-order inhomogeneity.
Thanks to our preparations provided by Corollary \ref{cor4} and Lemma \ref{lem6}, the resulting influence can appropriately
be estimated in terms of the expressions in (\ref{M}), and of the quantities under consideration themselves, namely of
\be{N}
	N_{q,\beta}(t):= 1+\sup_{s\in (0,t)} \Big\| B^\beta \Big(z(\cdot,s)-e^{-sB} z_0\Big) \Big\|_{L^q(\Om)},
	\qquad t\in (0,T),
\ee
at a conveniently controllable sublinear power (cf.~(\ref{7.18}).
If the fractional power $\beta$ appearing herein is suitably close to $1$, then this indeed yields the following key
result on regularity of $z$.
\begin{lem}\label{lem7}
  Let $\vartheta\in (\frac{3}{4},1)$, and suppose that $m>0, K>0, p\ge 2, q\ge 2$ and $\beta\in (\frac{1}{2},1)$ are such that
  (\ref{6.1}) holds, and that
  \be{7.1}
	2p(1-\beta) \ge 3.
  \ee
  Then for all $\eta>0$ there exists $C(m,K,p,q,\beta,\eta)>0$ with the property that if
  (\ref{f}), (\ref{init}), (\ref{vz}) and (\ref{0vz}) hold with some $T\in (0,\infty]$, then for all $t\in (0,T)$,
  \be{7.2}
	\Big\|B^\beta \Big(z(\cdot,t)-e^{-tB} z_0\Big) \Big\|_{L^q(\Om)}
	\le C(m,K,p,q,\beta,\eta) \cdot 
	\bigg\{ 1 + \sup_{s\in (0,t)} \|f(\cdot,s)\|_{L^p(\Om)} \bigg\}^{\frac{p}{p-1} \cdot (\frac{2q\beta+q-1}{3q} + \eta)}.
  \ee
\end{lem}
\proof
  We prepare our estimation procedure by firstly using that $q\ge \frac{3}{2\beta-1}$ and $p\ge \frac{3}{2(-\beta)}$ to verify that
  \be{7.22}
	\lam=\lam(p,q,\beta):=\min \Big\{ q \, , \, \frac{3q}{2q-2q\beta+1} \Big\}
  \ee
  satisfies
  \be{7.3}
	\lam>\min \Big\{ q \, , \, \frac{3q}{2q-2q\beta+3} \Big\}
	= \frac{3q}{2q-2q\beta+3}
  \ee
  and hence, in particular,
  \bas
	\lam>\frac{3}{2-2\beta+\frac{3}{q}}
	> \frac{3}{2-2\beta+(2\beta-1)}
	= 3>1,
  \eas
  and moreover has the property that
  \be{7.33}
	\lam< \frac{3q}{2q-2q\beta} = \frac{3}{2(1-\beta)} \le p,
  \ee
  so that
  \be{7.4}
	\|\varphi\|_{L^\lam(\Om)} \le \|\varphi\|_{L^p(\Om)}^\frac{p(\lam-1)}{(p-1)\lam} \|\varphi\|_{L^1(\Om)}^\frac{p-\lam}{(p-1)\lam}
	\qquad \mbox{for all } \varphi\in L^p(\Om)
  \ee
  by the H\"older inequality.\abs
  We next observe that, again by (\ref{6.1}), $\eta_0:=2q\beta-q-3$ is positive and thus
  \be{7.5}
	\psi(\wh{\eta}) := \bigg\{
	\frac{2q\beta-q-3+\wh{\eta}}{3q} + \wh{\eta} + \frac{(q-1)(2q\beta-q-3)}{6q^2\beta} \bigg\}
		\cdot \frac{2q\beta}{2q\beta-q-3-\wh{\eta}},
	\qquad \wh{\eta} \in (0,\eta_0),
  \ee
  well-defined with
  \bas
	\psi(\wh{\eta}) 
	&\searrow &
	\bigg\{ \frac{2q\beta-q-3}{3q} + \frac{(q-1)(2q\beta-q-3)}{6q^2\beta} \bigg\} \cdot \frac{2q\beta}{2q\beta-q-3} \\
	&=& \frac{2\beta}{3} + \frac{q-1}{3q} 
	= \frac{2q\beta+q-1}{3q}
	\qquad \mbox{as } \wh{\eta} \searrow 0,
  \eas
  whence given $\eta>0$ we can pick $\eta_1=\eta_1(p,q,\beta,\eta)\in (0,\eta_0)$ such that
  \be{7.6}
	\psi(\eta_1) \le \frac{2q\beta+q-1}{3q} + \eta,
  \ee
  where we can clearly achieve that, simultaneously,
  \be{7.7}
	\eta_1 < 2q-2q\beta.
  \ee
  According to the latter, namely, setting
  \be{7.8}
	r=r(p,q,\beta,\eta):=\frac{3q}{2q-2q\beta+3-\eta_1}
  \ee
  introduces a well-defined positive number $r$ which satisfies
  \be{7.9}
	r<\frac{3q}{2q-2q\beta+3-(2q-2q\beta)} =q,
  \ee
  and for which we moreover have
  \be{7.10}
	r>\frac{3q}{2q-2q\beta+3}>3,
  \ee
  again thanks to (\ref{6.1}).\abs
  We now recall known smoothing properties of the semigroup $(e^{-tB})_{t\ge 0}$ (\cite{friedman}, \cite{win_JDE}) to, firstly, fix
  $C_1>0$ such that
  \be{7.100}
	\|\na e^{-tB}\varphi\|_{W^{1,\infty}(\Om)} \le C_1 \|\varphi\|_{W^{1,\infty}(\Om)}
	\qquad \mbox{for all } \varphi\in W^{1,\infty}(\Om),
  \ee
  to, secondly, see that since $1<\lam\le q$ by (\ref{7.33}) and (\ref{7.22}), there exists $C_2=C_2(p,q,\beta)>0$ fulfilling
  \be{7.11}
	\|B^\beta e^{-tB} \varphi\|_{L^q(\Om)}
	\le C_2 \cdot (1+t^{-\beta-\frac{3}{2}(\frac{1}{\lam}-\frac{1}{q})}) e^{-t} \|\varphi\|_{L^\lam(\Om)}
	\qquad \mbox{for all $\varphi\in C^0(\bom)$ and } t>0,
  \ee
  and to, thirdly, observe that as $1<r<q$ by (\ref{7.10}) and (\ref{7.9}), we can similarly find $C_3=C_3(p,q,\beta,\eta)>0$ such that
  \be{7.12}
	\|B^\beta e^{-tB} \varphi\|_{L^q(\Om)}
	\le C_3 \cdot (1+t^{-\beta-\frac{3}{2}(\frac{1}{r}-\frac{1}{q})}) e^{-t} \|\varphi\|_{L^r(\Om)}
	\qquad \mbox{for all $\varphi\in C^0(\bom)$ and } t>0.
  \ee
  Apart from that, based on the fact that (\ref{7.10}) actually even warrants that $r>3$, we may employ Corollary \ref{cor4}
  and furthermore invoke Lemma \ref{lem6} to infer that whenever $m>0$ and $K>0$, we can find $C_i=C_i(m,K,p,q,\beta,\eta)>0$,
  $i\in\{4,5\}$, such that given $T\in (0,\infty]$ as well as $f,z_0,v_0,z,v$ and $Q$ satisfying 
  (\ref{f}), (\ref{init}), (\ref{vz}) and (\ref{0vz}), we have
  \be{7.13}
	\|v(\cdot,t)\|_{L^r(\Om)}
	\le C_4 M_p^{\frac{p}{p-1} \cdot (\frac{r-3}{3r} + \eta_1)}(t)
	\qquad \mbox{for all } t\in (0,T)
  \ee
  and
  \be{7.14}
	\Big\| \na \Big(z(\cdot,t)-e^{-tB} z_0\Big) \Big\|_{L^\infty(\Om)}
	\le C_5 M_p^{\frac{p}{p-1} \cdot \frac{(q-1)(2q\beta-q-3)}{6q^2\beta}}(t)
	\cdot N_{q,\beta}^\frac{q+3+\eta_1}{2q\beta}(t)
	\qquad \mbox{for all } t\in (0,T),
  \ee
  with $(M_p(t))_{t\in (0,T)}$ and $(N_{q,\beta}(t))_{t\in (0,T)}$ as given by (\ref{M}) and (\ref{N}).\abs
  We henceforth let $T\in (0,\infty]$ as well as $f,z_0,v_0,z,v$ and $Q$ be given such that 
  (\ref{f}), (\ref{init}), (\ref{vz}) and (\ref{0vz}) hold with some $m>0$ and $K>0$, and accordingly
  define $(M_p(t))_{t\in (0,T)}$ and $(N_{q,\beta}(t))_{t\in (0,T)}$ through (\ref{M}) and (\ref{N}).
  Then relying on a Duhamel representation associated with the first sub-problem in (\ref{0vz}), we can use (\ref{7.11}) and
  (\ref{7.12}) to estimate
  \bea{7.16}
	& & \hs{-20mm}
	\Big\| B^\beta \Big( z(\cdot,t)-e^{-tB}z_0\Big) \Big\|_{L^q(\Om)} \nn\\
	&=& \bigg\|
	\int_0^t B^\beta e^{-(t-s)B} f(\cdot,s) ds 
	- \int_0^t B^\beta e^{-(t-s)B} \Big\{ v(\cdot,t)\cdot\na z(\cdot,s)\Big\} ds \bigg\|_{L^q(\Om)} \nn\\
	&\le& C_2 \int_0^t \Big(1+(t-s)^{-\beta-\frac{3}{2}(\frac{1}{\lam}-\frac{1}{q})} \Big) 
		e^{-(t-s)} \|f(\cdot,s)\|_{L^\lam(\Om)} ds \nn\\
	& & + C_3 \int_0^t \Big(1+(t-s)^{-\beta-\frac{3}{2}(\frac{1}{1}-\frac{1}{q})} \Big) e^{-(t-s)} 
		\|v(\cdot,s)\cdot\na z(\cdot,s)\|_{L^r(\Om)} ds
  \eea
  for $t\in (0,T)$, where by (\ref{7.4}) and (\ref{f}),
  \be{7.17}
	\|f(\cdot,t)\|_{L^\lam(\Om)}
	\le m^\frac{p-\lam}{(p-1)\lam} M_p^\frac{p(\lam-1)}{(p-1)\lam}(t)
	\qquad \mbox{for all $t\in (0,T)$ and } s\in (0,t),
  \ee
  and where thanks to (\ref{7.13}), (\ref{7.14}), (\ref{7.100}), (\ref{init}) and the fact that $M_p\ge 1$ and $N_{q,\beta}\ge 1$,
  for all $t\in (0,T)$ and $s\in (0,t)$ we have
  \bea{7.18}
	\|v(\cdot,s)\cdot\na z(\cdot,s)\|_{L^r(\Om)}
	&\le& \|v(\cdot,s)\|_{L^r(\Om)} \|\na z(\cdot,s)\|_{L^\infty(\Om)} \nn\\
	&\le& \|v(\cdot,s)\|_{L^r(\Om)} \cdot \bigg\{
	\Big\| \na \Big(z(\cdot,s)-e^{-sB} z_0\Big) \Big\|_{L^\infty(\Om)}
	+ \|\na e^{-sB} z_0\|_{L^\infty(\Om)} \bigg\} \nn\\
	&\le& C_4 M_p^{\frac{p}{p-1} \cdot (\frac{r-3}{3r} + \eta_1)}(t)
	\cdot \Big\{ C_5 M_p^{\frac{p}{p-1} \cdot \frac{(q-1)(2q\beta-q-3)}{6q^2\beta}}(t) \cdot 
		N_{q,\beta}^\frac{q+3+\eta_1}{2q\beta}(t) 
	+ C_1 K \Big\} \nn\\
	&\le& C_6 M_p^{\frac{p}{p-1} \cdot (\frac{r-3}{3r} + \eta_1 + \frac{(q-1)(2q\beta-q-3)}{6q^2\beta})}(t)
		\cdot N_{q,\beta}^\frac{q+3+\eta_1}{2q\beta}(t)
  \eea
  with $C_6=C_6(m,K,p,q,\beta,\eta):=C_4 \cdot \max \{C_5 \, , \, C_1 K\}$.
  Since
  \bas
	C_7=C_7(p,q,\beta):=\int_0^\infty (1+\sigma^{-\beta-\frac{3}{2}(\frac{1}{\lam}-\frac{1}{q})}) e^{-\sigma} d\sigma
	\quad \mbox{and} \quad
	C_8=C_8(p,q,\beta,\eta):=\int_0^\infty (1+\sigma^{-\beta-\frac{3}{2}(\frac{1}{r}-\frac{1}{q})}) e^{-\sigma} d\sigma
  \eas
  are both finite due to the circumstance that
  \bas
	\beta+\frac{3}{2}\Big(\frac{1}{\lam}-\frac{1}{q}\Big)
	< \beta + \frac{3}{2} \cdot \Big(\frac{2q-2q\beta+3}{3q}-\frac{1}{q}\Big)=1
  \eas
  by (\ref{7.22}), and thet, similarly, $\beta+\frac{3}{2}(\frac{1}{r}-\frac{1}{q})<1$ by the first inequality in (\ref{7.10}),
  from (\ref{7.16})-(\ref{7.18}) we thus conclude that for all $t\in (0,T)$,
  \be{7.19}
	\Big\| B^\beta \Big( z(\cdot,t)-e^{-tB}z_0\Big) \Big\|_{L^q(\Om)} 
	\le C_9 M_p^\frac{p(\lam-1)}{(p-1)\lam}(t)
	+ C_9 M_p^{\frac{p}{p-1} \cdot (\frac{r-3}{3r} + \eta_1 + \frac{(q-1)(2q\beta-q-3)}{6q^2\beta})}(t)
		\cdot N_{q,\beta}^\frac{q+3+\eta_1}{2q\beta}(t)
  \ee
  with
  \bas
	C_9=C_9(m,K,p,q,\beta,\eta):=\max \big\{ C_2 C_7 m^\frac{p-\lam}{(p-1)\lam} \, , \, C_3 C_6 C_8\big\}.
  \eas
  We may now rely on the inequality 
  $\frac{q+3+\eta_1}{2q\beta} < \frac{q+3+(2q\beta-q-3)}{2q\beta}=1$, as asserted by our restriction that $\eta_1<\eta_0$,
  to see that due to Young's inequality there exists $C_{10}=C_{10}(m,K,p,q,\beta,\eta)>0$ such that
  \bas
	ab \le \frac{1}{2} a^\frac{2q\beta}{q+3+\eta_1} + C_{10} b^\frac{2q\beta}{2q\beta-q-3-\eta_1}
	\qquad \mbox{for all $a\ge 0$ and } b\ge 0.
  \eas
  From (\ref{7.19}) we therefore obtain that for all $t\in (0,T)$,
  \bas
	N_{q,\beta}(t)
	&\le& 1+ C_9 M_p^\frac{p(\lam-1)}{(p-1)\lam}(t) \\
	& & + \frac{1}{2} N_{q,\beta}(t)
	+ C_{10} \cdot 
	\Big\{ C_9 M_p^{\frac{p}{p-1} \cdot (\frac{r-3}{3r} + \eta_1 + \frac{(q-1)(2q\beta-q-3)}{6q^2\beta})}(t)
	\Big\}^\frac{2q\beta}{2q\beta-q-3-\eta_1},
  \eas
  and that thus
  \bea{7.20}
	\hs{-8mm}
	N_{q,\beta}(t)
	&\le& 2+ 2 C_9 M_p^\frac{p(\lam-1)}{(p-1)\lam}(t) \\
	& & 
	+ 2C_9^\frac{2q\beta}{2q\beta-q-3-\eta_1} C_{10} 
	M_p^{\frac{p}{p-1} \cdot (\frac{r-3}{3r} + \eta_1 + \frac{(q-1)(2q\beta-q-3)}{6q^2\beta}) 
		\cdot \frac{2q\beta}{2q\beta-q-3-\eta_1}}(t)
	\qquad \mbox{for all } t\in (0,T).
  \eea
  Since
  \bas
	\frac{\lam-1}{\lam}
	= 1 - \frac{1}{\lam}
	\le 1 - \frac{2q-2q\beta+1}{3q}
	= \frac{2q\beta+q-1}{3q} 
	\le \frac{2q\beta+q-1}{3q} +\eta
  \eas
  thanks to (\ref{7.22}), and since (\ref{7.8}) says that
  \bas
	\frac{r-3}{3r}
	= \frac{1}{3} - \frac{1}{r} 
	= \frac{1}{3} - \frac{2q-2q\beta+3-\eta_1}{3q}
	= \frac{2q\beta-q-3+\eta_1}{3q},
  \eas
  and that therefore
  \bas
	\Big( \frac{r-3}{3r} + \eta_1 + \frac{(q-1)(2q\beta-q-3)}{6q^2\beta} \Big) \cdot \frac{2q\beta}{2q\beta-q-3-\eta_1}
	= \psi(\eta_1) \le \frac{2q\beta+q-1}{3q} + \eta
  \eas
  according to (\ref{7.6}), according to the definitions of $(M_p(t))_{t\in (0,T)}$ and $(N_{q,\beta}(t))_{t\in (0,T)}$
  in (\ref{M}) and (\ref{N}) we directly infer (\ref{7.2}) from (\ref{7.20}).
\qed
Our main result on (\ref{0vz}) can now be achieved by suitably adjusting the auxiliarly parameters $q$ and $\beta$ in the above,
and by adequately coping with the correction term $e^{-tB} z_0$ in (\ref{7.2}):\abs
\proofc of Theorem \ref{theo2}. \quad
  Since $\frac{2p-3}{2p}=1-\frac{3}{2p}>\frac{1}{2}$ due to our assumption that $p>3$, we can fix $\beta=\beta(p)>\frac{1}{2}$ such that
  $\beta\le \frac{2p-3}{2p}$, and that thus (\ref{7.1}) holds.
  Given $\eta>0$, we thereupon pick $q=q(p)>\frac{3}{2\beta-1}$ suitably large fulfilling
  \bas
	\frac{2(q+1)}{3q}<\frac{2}{3}+\eta,
  \eas
  which ensures that
  \be{2.3}
	\psi(\wh{\eta}):=
	\frac{(q-1)(2q\beta-q-3)}{6q^2\beta} + \Big( \frac{2q\beta+q-1}{3q} + \wh{\eta}\Big) \cdot 
	\Big( \frac{q+3}{2q\beta}+\wh{\eta}\Big),
	\qquad \wh{\eta}>0,
  \ee
  satisfies
  \bas
	\psi(\wh{\eta}) \searrow 
	\frac{(q-1)(2q\beta-q-3)}{6q^2\beta} + \frac{(q+3)(2q\beta+q-1)}{6q^2\beta}
	=\frac{2(q+1)}{3q}<\frac{2}{3}+\eta,
  \eas
  and that hence we can choose $\eta_1=\eta_1(p,\eta)>0$ in such a way that
  \be{2.4}
	\psi(\eta_1) \le \frac{2}{3} + \eta.
  \ee
  We now employ Lemma \ref{lem6} and Lemma \ref{lem7} to see that given $m>0$ and $K>0$ we can find $C_1=C_1(m,K,p,\eta)>0$
  such that whenever $T\in (0,\infty]$ and $f$ as well as $z_0,v_0,z,v$ and $Q$ comply with 
  (\ref{f}), (\ref{init}), (\ref{vz}) and (\ref{0vz}),
  as before letting $(M_p(t))_{t\in (0,T)}$ and $N_{q,\beta}(t))_{t\in (0,T)}$ be as defined in (\ref{M}) and (\ref{N}) we have
  \be{2.5}
	\|z(\cdot,t)-e^{-tB}z_0\|_{W^{1,\infty}(\Om)}
	\le C_1 M_p^{\frac{p}{p-1} \cdot \frac{(q-1)(2q\beta-q-3)}{6q^2\beta}}(t) 
		\cdot N_{q,\beta}^{\frac{q+3}{2q\beta}+\eta_1}(t)
	\qquad \mbox{for all } t\in (0,T)
  \ee
  and
  \be{2.6}
	N_{q,\beta}(t)
	\le C_2 M_p^{\frac{p}{p-1} \cdot (\frac{2q\beta+q-1}{3q} + \eta_1)}(t)
	\qquad \mbox{for all } t\in (0,T).
  \ee
  Finally fixing $C_3>0$ such that, in accordance with a known boundedness feature of the Neumann heat semigroup,
  \bas
	\|e^{-tB}\varphi\|_{W^{1,\infty}(\Om)} 
	\le C_3 \|\varphi\|_{W^{1,\infty}(\Om)}
	\qquad \mbox{for all } \varphi\in W^{1,\infty}(\Om),
  \eas
  we may combine (\ref{2.5}) with (\ref{2.6}) and (\ref{init}) to infer that given any such $m,K,T,f,z_0,v_0,z$ and $v$,
  and correspondingly taking $(M_p(t))_{t\in (0,T)}$ and $N_{q,\beta}(t))_{t\in (0,T)}$ from (\ref{M}) and (\ref{N}), we can estimate
  \bas
	\|z(\cdot,t)\|_{W^{1,\infty}(\Om)}
	&\le& \|z(\cdot,t)-e^{-tB}z_0\|_{W^{1,\infty}(\Om)}
	+ \|e^{-tB} z_0\|_{W^{1,\infty}(\Om)} \\
	&\le& C_1 M_p^{\frac{p}{p-1} \cdot \frac{(q-1)(2q\beta-q-3)}{6q^2\beta}}(t) \cdot N_{q,\beta}^{\frac{q+3}{2q\beta}+\eta_1)}(t)
	+ C_3 \|z_0\|_{W^{1,\infty}(\Om)} \\
	&\le& C_1 C_2^{\frac{q+3}{2q\beta}+\eta_1} 
	M_p^{\frac{p}{p-1} \cdot \big\{ \frac{(q-1)(2q\beta-q-3)}{6q^2\beta} + (\frac{2q\beta+q-1}{3q}+\eta_1) \cdot 
		(\frac{q+3}{2q\beta}+\eta_1)\big\}}(t)
	+ C_3 K \\
	&\le& C_1 C_2^{\frac{q+3}{2q\beta}+\eta_1} 
	M_p^{\frac{p}{p-1}\cdot \psi(\eta_1)}(t)
	+ C_3 K
	\qquad \mbox{for all } t\in (0,T)
  \eas
  by (\ref{2.3}).
  In view of (\ref{2.4}), this immediately establishes (\ref{2.1}).
\qed
\mysection{Boundedness in (\ref{0}) for all subcritical nonlinearities. Proof of Theorem \ref{theo11}}\label{sect3}
Next addressing the flux-limited Keller-Segel-Stokes system (\ref{0}), let us first recall standard theory on 
problems of related types to state the following basic result on local existence and extensibility, actually available
for fairly general $F$.
\begin{lem}\label{lem12}
  If $\Phi\in W^{2,\infty}(\Om)$ and $F\in C^2([0,\infty))$,
  and if $n_0,c_0$ and $u_0$ comply with (\ref{incu}), then there exist 
  $\tm\in (0,\infty]$ and unique functions
  \bas
	\left\{ \begin{array}{l}
	n\in C^0(\bom\times [0,\tm)) \cap C^{2,1}(\bom\times (0,\tm)), \\[1mm]
	c\in \bigcap_{q>3} C^0([0,\tm);W^{1,q}(\Omega)) \cap C^{2,1}(\bom\times (0,\tm))
	\qquad \mbox{and} \\[1mm]
	u\in \bigcup_{\vartheta\in (\frac{1}{2},1)} C^0([0,\tm);D(A^\vartheta)) \cap C^{2,1}(\bom\times (0,\tm);\R^3)
	\end{array} \right.
  \eas
  with the properties that $n\ge 0$ and $c\ge 0$ in $\Om\times (0,\tm)$, that one can find
  $P\in C^{1,0}(\Om\times (0,\tm))$ such that $(n,c,u,P)$ is a classical solution of (\ref{0}) in $\Om\times (0,\tm)$, and that
  \bea{ext}
	& & \hs{-20mm}
	\mbox{if $\tm=\infty$, \quad then for all $\vartheta\in (\frac{3}{4},1)$,} \nn\\
	& & \limsup_{t\nearrow \tm} \Big\{ 
		\|n(\cdot,t)\|_{L^\infty(\Om)} + \|c(\cdot,t)\|_{W^{1,\infty}(\Om)} + \|A^\vartheta u(\cdot,t)\|_{L^2(\Om)}
		\Big\} = \infty.
  \eea
  For this solution we additionally have
  \be{mass}
	\io n(\cdot,t)=\io n_0
	\qquad \mbox{for all } t\in (0,\tm).
  \ee
\end{lem}
\proof
  This can be seen by standard arguments well-documented in closely related contexts (cf.~\cite{win_CPDE}, for instance).
\qed
Now the reward for all our efforts related to Theorem \ref{theo2} consists in the circumstance that
its outcome facilitates the derivation of $L^p$ bounds for the first solution component in (\ref{0}) through
a noticeably short argument which in essence reduces to quite a straighforward combination of a standard
testing procedure with the conditional estimate in (\ref{2.1}).
\begin{lem}\label{lem13}
  Suppose that $\Phi\in W^{2,\infty}(\Om)$, that $F\in C^2([0,\infty))$ satisfies (\ref{F}) with some $\Kf>0$ and $\alpha>\frac{1}{2}$,
  and that (\ref{incu}) holds.
  Then the solution of (\ref{0}) from Lemma \ref{lem12} has the property that
  \be{13.1}
	\sup_{t\in (0,\tm)} \|n(\cdot,t)\|_{L^p(\Om)} <\infty
	\qquad \mbox{for all } p>3.
  \ee
\end{lem}
\proof
  Without loss of generality assuming that $\alpha<1$, we note that since $2\cdot (1-\alpha)<1$ by hypothesis, we can choose
  $\eta>0$ such that still
  \be{13.2}
	(2+3\eta)\cdot (1-\alpha) <1.
  \ee
  Taking $\vartheta\in (\frac{3}{4},1)$ such that $u_0\in D(A^\vartheta)$, 
  for fixed $p>3$ we may then apply Theorem \ref{theo2} to $m:=1+\io n_0>0$, $K:=1+\max\{\|c_0\|_{W^{1,\infty}(\Om)} \, , \,
  \|A^\vartheta u_0\|_{L^2(\Om)} \}>0$, $(f,z,v,Q):=(n,c,u,P)$ and $T:=\tm$ to find $C_1=C_1(p)>0$ such that
  \be{13.3}
	\|\na c(\cdot,t)\|_{L^\infty(\Om)}
	\le C_1 M_p^{\frac{p}{p-1} \cdot (\frac{2}{3}+\eta)}(t)
	\qquad \mbox{for all } t\in (0,\tm),
  \ee
  where in line with (\ref{M}), we have set $M_p(t):=1+\sup_{s\in (0,t)} \|n(\cdot,s)\|_{L^p(\Om)}$ for $t\in (0,\tm)$.\abs
  To derive (\ref{13.1}) from this, we test the first equation in (\ref{0}) against $n^{p-1}$ and intergate by parts to obtain
  that since $\na\cdot u=0$, due to Young's inequality we have
  \bas
	\frac{1}{p} \frac{d}{dt} \io n^p + (p-1) \io n^{p-2} |\na n|^2 	
	&=& (p-1) \io n^{p-1} F(|\na c|^2) \na n\cdot\na c \\
	&\le& \frac{p-1}{2} \io n^{p-2}|\na n|^2
	+ \frac{p-1}{2} \io n^p F^2(|\na c|^2) |\na c|^2 \\
	&\le& \frac{p-1}{2} \io n^{p-2}|\na n|^2
	+ \frac{(p-1)\Kf^2}{2} \io n^p |\na c|^{2-2\alpha}
  \eas
  for all $t\in (0,\tm)$, so that
  \be{13.4}
	\frac{d}{dt} \io n^p
	+ \frac{2(p-1)}{p} \io |\na n^\frac{p}{2}|^2
	\le \frac{p(p-1)\Kf^2}{2} \io n^p |\na c|^{2-2\alpha}
	\qquad \mbox{for all } t\in (0,\tm).
  \ee
  Here since we are assuming $2-2\alpha$ to be positive, we may utilize (\ref{13.3}) to estimate
  \bea{13.5}
	\frac{p(p-1)\Kf^2}{2} \io n^p |\na c|^{2-2\alpha}
	&\le& \frac{p(p-1)\Kf^2}{2} \|\na c\|_{L^\infty(\Om)}^{2-2\alpha} \io n^p \nn\\
	&\le& C_2 M_p^{\frac{p}{p-1} \cdot (\frac{2}{3}+\eta)\cdot (1-2\alpha)}(t) \io n^p
	\qquad \mbox{for all } t\in (0,\tm),
  \eea
  where $C_2=C_2(p):=\frac{p(p-1)\Kf^2}{2} \cdot C_1^{2-2\alpha}$.
  Now thanks to (\ref{mass}), an interpolation on the basis of the Gagliardo-Nirenberg inequality shows that with some
  $C_3=C_3(p)>0$ and $C_4=C_4(p)>0$ we have
  \bas
	\bigg\{ \io n^p \bigg\}^\frac{3p-1}{3(p-1)}
	&=& \|n^\frac{p}{2}\|_{L^2(\Om)}^\frac{2(3p-1)}{3(p-1)} \\
	&\le& C_3 \|\na n^\frac{p}{2}\|_{L^2(\Om)}^2 \|n^\frac{p}{2}\|_{L^\frac{2}{p}(\Om)}^\frac{4}{3(p-1)}
	+ C_3 \|n^\frac{p}{2}\|_{L^\frac{2}{p}(\Om)}^\frac{2(3p-1)}{3(p-1)} \\
	&\le& C_4 \|\na n^\frac{p}{2}\|_{L^2(\Om)}^2 + C_4
	\qquad \mbox{for all } t\in (0,\tm),
  \eas
  and that thus
  \be{13.6}
	\frac{2(p-1)}{p} \io |\na n^\frac{p}{2}|^2
	\ge C_5 \cdot \bigg\{ \io n^p \bigg\}^\kappa - C_5
	\qquad \mbox{for all } t\in (0,\tm)
  \ee
  with $C_5=C_5(p):=\frac{2(p-1)}{pC_4}$ and 
  \be{13.66}
	\kappa=\kappa(p):=\frac{3p-1}{3(p-1)}>1.
  \ee
  Now a combination of (\ref{13.6}) with (\ref{13.5}) and Young's inequality shows that for each $T\in (0,\tm)$, (\ref{13.4}) entails
  the inequality
  \bas
	\frac{d}{dt} \io n^p
	+ C_5 \cdot \bigg\{ \io n^p \bigg\}^\kappa
	&\le& C_2 M_p^{\frac{p}{p-1} \cdot (\frac{2}{3}+\eta)\cdot (2-2\alpha)}(T) \io n^p + C_5 \nn\\
	&=& \Bigg\{ \frac{C_5}{2} \cdot \bigg\{ \io n^p \bigg\}^\kappa \Bigg\}^\frac{1}{\kappa} \cdot
	\Bigg\{ \Big(\frac{2}{C_5}\Big)^\frac{1}{\kappa} \cdot C_2 M_p^{\frac{p}{p-1}\cdot (\frac{2}{3}+\eta)\cdot (2-2\alpha)}(T)
		\Bigg\}	
	+ C_5 \nn\\
	&\le& \frac{C_5}{2} \cdot \bigg\{ \io n^p \bigg\}^\kappa \nn\\
	& & + C_6 M_p^{\frac{\kappa}{\kappa-1} \cdot \frac{p}{p-1}\cdot (\frac{2}{3}+\eta)\cdot (2-2\alpha)}(T) + C_5
	\qquad \mbox{for all } t\in (0,T)
  \eas
  with $C_6=C_6(p):=\big\{ (\frac{2}{C_5})^\frac{1}{\kappa} C_2\big\}^\frac{\kappa}{\kappa-1}$.
  Since $M_p\ge 1$, this implies that if we let $C_7=C_7(p):=C_6+C_5$, then for any such $T$ we have
  \bas
	\frac{d}{dt} \io n^p + \frac{C_5}{2} \cdot \bigg\{ \io n^p \bigg\}^\kappa
	\le C_7  M_p^{\frac{\kappa}{\kappa-1} \cdot \frac{p}{p-1}\cdot (\frac{2}{3}+\eta)\cdot (2-2\alpha)}(T)
	\qquad \mbox{for all } t\in (0,T),
  \eas
  which through an ODE comparison argument guarantees that
  \bas
	\io n^p(\cdot,t)
	\le \max \Bigg\{ \io n_0^p \, , \, 
	\bigg\{ \frac{2C_7}{C_5}  M_p^{\frac{\kappa}{\kappa-1} \cdot \frac{p}{p-1}\cdot (\frac{2}{3}+\eta)\cdot (2-2\alpha)}(T)
		\bigg\}^\frac{1}{\kappa} \Bigg\}
	\qquad \mbox{for all } t\in (0,T),
  \eas
  and that therefore, again since $M_p\ge 1$,
  \bea{13.8}
	M_p(T)
	&\le&
	1 + \max \bigg\{ \|n_0\|_{L^p(\Om)} \, , \, \Big(\frac{2C_7}{C_5}\Big)^\frac{1}{p\kappa} M_p^\lam(T) \bigg\} \nn\\
	&\le& C_8 M_p^\lam(T)
	\qquad \mbox{for all } T\in (0,\tm)
  \eea
  with $C_8=C_8(p):=1+\max \big\{ \|n_0\|_{L^p(\Om)} \, , \, (\frac{2C_7}{C_5})^\frac{1}{p\kappa}\big\}$
  and $\lam=\lam(p):=\frac{1}{\kappa-1} \cdot \frac{1}{p-1} \cdot (\frac{2}{3}+\eta)\cdot (2-2\alpha)$.\abs
  It now only remains to observe that according to (\ref{13.66}), our restriction on $\eta$ in (\ref{13.2}) ensures that
  \bas
	\lam
	= \frac{3(p-1)}{2} \cdot \frac{1}{p-1} \cdot \Big(\frac{2}{3}+\eta\Big)\cdot (2-2\alpha) 
	= (2+3\eta)\cdot (1-\alpha)
	< 1
  \eas
  to finally conclude from (\ref{13.8}) that
  \bas
	M_p(T) \le C_8^\frac{1}{1-\lam}
	\qquad \mbox{for all } T\in (0,\tm),
  \eas
  which implies (\ref{13.1}) upon taking $T\nearrow \tm$.
\qed
In view of Theorem \ref{theo2} and Proposition \ref{prop1}, the latter immediately implies bounds for the 
quantities in (\ref{ext}) related to the signal concentration and the fluid velocity.
\begin{lem}\label{lem14}
  If $\Phi\in W^{2,\infty}(\Om)$ and $F\in C^2([0,\infty))$ is such that (\ref{F}) is valid with some $\Kf>0$ and $\alpha>\frac{1}{2}$,
  and if (\ref{incu}) holds, 
  then there exists $\vartheta>\frac{3}{4}$ such that $\tm$ as well as the functions $c$ and $u$ from Lemma \ref{lem12} satisfy
  \be{14.1}
	\sup_{t\in (0,\tm)} \|c(\cdot,t)\|_{W^{1,\infty}(\Om)} <\infty
  \ee
  and
  \be{14.2}
	\sup_{t\in (0,\tm)} \|A^\vartheta u(\cdot,t)\|_{L^2(\Om)} <\infty.
  \ee
\end{lem}
\proof
  The boundedness feature in (\ref{14.1}) is an immediate consequence of Theorem \ref{theo2} when applied to any fixed $p>3$ and 
  combined with Lemma \ref{lem13}, while (\ref{14.2}) similarly results from Proposition \ref{prop1}.
\qed
Now thanks to the $L^\infty$ bounds for $\na c$ and $u$ implied by Lemma \ref{lem14} due to the continuity of the embedding
$D(A^\vartheta) \hra L^\infty(\Om;\R^3)$, a straightforward applciation of heat semigroup estimates to the first
equation in (\ref{0}) finally yields $L^\infty$ estimates also for $n$:
\begin{lem}\label{lem15}
  Let $\Phi\in W^{2,\infty}(\Om)$ and $F\in C^2([0,\infty))$ be such that (\ref{F}) holds with some $\Kf>0$ and $\alpha>\frac{1}{2}$,
  and assume (\ref{incu}).
  Then with $\tm$ and $n$ taken from Lemma \ref{lem12}, we have
  \be{15.1}
	\sup_{t\in (0,\tm)} \|n(\cdot,t)\|_{L^\infty(\Om)} <\infty.
  \ee
\end{lem}
\proof
  We fix any $\lam>3$ and then readily infer from (\ref{F}), Lemma \ref{lem13} and Lemma \ref{lem14} the existence of
  $C_1>0$ and $C_2>0$ such that $h_1:=nF(|\na c|^2)\na c + nu$ and $h_2:=n$ satisfy
  \bas
	\|h_1(\cdot,t)\|_{L^\lam(\Om)} \le C_1
	\quad \mbox{and} \quad
	\|h_2(\cdot,t)\|_{L^\frac{\lam}{2}(\Om)} \le C_2
	\qquad \mbox{for all } t\in (0,\tm).
  \eas
  As known smoothing estimates for the Neumann heat semigroup $(e^{t\Delta})_{t\ge 0}$ on $\Om$ (\cite{FIWY}, \cite{win_JDE})
  provide $C_3>0$ and $C_4>0$ fulfilling
  \bas
	\|e^{t\Delta}\na\cdot\varphi\|_{L^\infty(\Om)}
	\le C_3 \cdot (1+t^{-\frac{1}{2}-\frac{3}{2\lam}}) \|\varphi\|_{L^\lam(\Omega)}
	\qquad \mbox{for all $\varphi\in C^1(\bom;\R^3)$ such that $\varphi\cdot\nu=0$ on } \pO
  \eas
  as well as
  \bas
	\|e^{t\Delta}\varphi\|_{L^\infty(\Om)}
	\le C_4 \cdot (1+t^{-\frac{3}{\lam}}) \|\varphi\|_{L^\frac{\lam}{2}(\Om)} 
	\qquad \mbox{for all } \varphi\in C^0(\bom),
  \eas
  by means of a variation-of-constants representation associated with the identity $n_t = \Delta n-n-\na\cdot h_1 + h_2$ we thus
  obtain that due to the maximum principle,
  \bas
	\|n(\cdot,t)\|_{L^\infty(\Om)}
	&=& \bigg\| e^{t(\Delta-1)} n_0
	- \int_0^t e^{(t-s)(\Delta-1)} \na \cdot h_1(\cdot,s) ds 
	+ \int_0^t e^{(t-s)(\Delta-1)} h_2(\cdot,s) ds \bigg\|_{L^\infty(\Om)} \\
	&\le& e^{-t} \|n_0\|_{L^\infty(\Om)}
	+ C_3 \int_0^t \Big( 1+(t-s)^{-\frac{1}{2}-\frac{3}{2\lam}} \Big) e^{-(t-s)} \|h_1(\cdot,s)\|_{L^\lam(\Om)} ds \\
	& & + C_4 \int_0^t \Big( 1+(t-s)^{-\frac{3}{\lam}} \Big) e^{-(t-s)} \|h_2(\cdot,s)\|_{L^\frac{\lam}{2}(\Om)} ds \\
	&\le& \|n_0\|_{L^\infty(\Om)}
	+ C_1 C_3 \int_0^\infty (1+\sigma^{-\frac{1}{2}-\frac{3}{2\lam}}) e^{-\sigma} d\sigma \\
	& & + C_2 C_4 \int_0^\infty (1+\sigma^{-\frac{3}{\lam}}) e^{-\sigma} d\sigma 
	\qquad \mbox{for all } t\in (0,\tm).
  \eas
  The claim thus follows from the observation that the rightmost two integrals herein are both finite due to the fact that
  the inequality $\lam>3$ warrants that $-\frac{1}{2}-\frac{3}{2\lam}>-1$ and $-\frac{3\{\lam}>-1$.
\qed
Our main result for (\ref{0}) has thus actually been established already:\abs
\proofc of Theorem \ref{theo11}.\quad
  Both the statement on global existence and uniqueness and the claim concerning the boundedness property in (\ref{11.1})
  directly result from Lemma \ref{lem12} in conjunction with Lemma \ref{lem14} and Lemma \ref{lem5}.
\qed
\vspace*{5mm}
{\bf Acknowledgement.} \quad
  The author acknowledges support of the {\em Deutsche Forschungsgemeinschaft} 
  in the context of the project {\em Emergence of structures and advantages in cross-diffusion systems} 
  (Project No.~411007140, GZ: WI 3707/5-1).

\end{document}